\documentclass[12pt]{article}
\usepackage{amsfonts}
\usepackage{amssymb}
\usepackage{textcomp}

\font\aaa=cmsy10 at 10pt
\renewcommand{\S}{{\aaa\char120\ }}

\renewcommand{\phi}{\varphi}
\newcommand{\be}{\begin{equation}}
\newcommand{\ee}{\end{equation}}
\newcommand{\ba}{\begin{eqnarray}}
\newcommand{\ea}{\end{eqnarray}}
\newcommand{\ban}{\begin{eqnarray*}}
\newcommand{\ean}{\end{eqnarray*}}

\newcommand{\too}{\mathop{\longrightarrow}}
\newcommand{\nul}{{\bf0}}
\newcommand{\odin}{{\bf1}}

\newcommand{\rd}{{\Bbb R}^d}

\newcommand{\zd}{{\Bbb Z}^{d}}

\renewcommand{\r}{{\Bbb R}}
\newcommand{\z} {{\Bbb Z}}
\newcommand{\n} {{\Bbb N}}

\newcommand{\ddd}{,\dots,}
\renewcommand{\lll}{\left(}
\newcommand{\rrr}{\right)}
\newcommand{\ex}[1]{e^{2\pi i{#1}}}
\newcommand{\exm}[1]{e^{-2\pi i{#1}}}
\newcommand{\sml}[3]{\sum\limits_{{#1}={#2}}^{#3}}
\newcommand{\Ldvad}{L_2(\rd)}
\newcommand{\h}{\widehat}

\newtheorem{theo}{Theorem}
\newtheorem{lem}[theo]{Lemma}
\newtheorem {prop} [theo] {Proposition}
\newtheorem {coro} [theo] {Corollary}
\newtheorem {defi} [theo] {Definition}

\begin {document}

\noindent

\vspace {.5cm}
\begin {center}
{ \LARGE\bf
Multivariate Wavelet Frames
\footnote {The paper is supported by RFBR, project N 06-01-00457

Unitary Extension Principle

Mathematics Subject Classifications 42C40}
}
\end {center}

\begin {center}
 \bf M. Skopina
\end {center}
\vspace {.5cm}

\begin {abstract}
We proved that for any
 matrix  dilation  and for any positive integer $n$,
 there exists a compactly supported tight  wavelet frame
  with  approximation order $n$.
Explicit  methods for construction of  dual and
tight wavelet frames with a given number of
vanishing moments are suggested.
\end {abstract}


\subsection {1. Introduction}

Approximation order  of a wavelet frame decomposition
\be
\sml
j{-\infty}\infty\sum\limits_{k\in\zd}\sml\nu1r \langle
f,\widetilde\psi_{jk}^{(\nu)}\rangle \psi_{jk}^{(\nu)}
\label{00}
\ee
is closely related to the number of   vanishing moments.
 It follows from Proposition 5.7 of~\cite{35} that a  MRA-based
compactly supported dual wavelet frame system with diagonal matrix dilation provides
approximation order $N$ whenever
the Fourier transform of all its generators 
has zero up to order $N$ at the origin (in our terminology, this means that
vanishing moment property of order $N-1$ holds for all the generators).
We give a proof of another fact: for a wide class of dual wavelet frames
( not necessary MRA-based, not necessary with diagonal matrix dilation
and not necessary compactly supported), vanishing moment property
for $\widetilde\psi^{(\nu)}, \nu=1\ddd r,$  is sufficient for the corresponding approximation
order of~(\ref{00}).  On the other hand, vanishing moment property of order $0$
for $\psi^{(\nu)}, \nu=1\ddd r,$ is a necessary condition for the system
 $\{\psi_{jk}^{(\nu)}\}$ to be a frame. This explains why  vanishing moments property is useful
 for frames,  and  why study of this topic is important.
 Some additional arguments are presented in Section~2.

 If a MRA-based wavelet system is a basis, then
the number of its vanishing moments depends only on the dual generating refinable function.
Situation is essentially different for frames. Two pairs of dual wavelet frames
may be generated by the same refinable functions and have different number of vanishing moments.
The goal of this paper is to describe  refinable functions generating dual wavelet
systems (potential frames)  with
vanishing moments and to  present  an explicit  method for construction of
compactly supported wavelet frames with arbitrary number of vanishing moments.
A close problem was studied by
Ming-Jun Lai and A. Petukhov\cite{MP} for univariate dual wavelet frames.
Their technique is not appropriate for multi-dimensional investigations because
 zero properties of  multivariate  masks can not be described by
means of factorization in contrast to the one-dimensional case.



Throughout the paper we will use the following notations.

$ \n $ is the set of positive integers,
$ \rd $ denotes the $d $-dimensional Euclidean space,
$x = (x_1\ddd x_d) $, $y = (y_1\ddd y_d) $ are its elements (vectors),
$ (x, y) = x_1y_1 +\dots+x_dy_d $,
$ |x | =\sqrt {(x, x)} $,
${\bf e}_j=(0\ddd1\ddd0)$ is the $j$-th
unit vector in $\rd$, \ $ \nul = (0\ddd0) \in \rd $, $ \odin = (1\ddd1) \in \rd $;
$ \zd $ is the integer lattice in $ \rd $.
For $x,y\in\rd$, we  write $x>y$ if $x_j>y_j$, $j=1\ddd d$;
$ \zd_+ = \{x\in Z^d:\ x\ge\nul \} $.
If $\alpha, \beta\in\zd_+$, $a,b\in\rd$, we set
$\alpha!=\prod\limits_{j=1}^d\alpha_j!$,
$\lll\alpha\atop\beta\rrr=
\frac{\alpha!}{\beta!(\alpha-\beta)!}$,
$a^b=\prod\limits_{j=1}^d{a_j}^{b_j}$,
 $[\alpha]=\sml j1d \alpha_j$,
$D^\alpha f=\frac{\partial^{[\alpha]}f}
{\partial^{\alpha_1}x_1\dots\partial^{\alpha_d}x_d}$;
$\delta_{ab}$ denotes  Kronecker delta; $\Bbb T^d$ is the unit $d$-dimensional torus;
$\Bbb C$ is the set of complex numbers.

Let $ M $ be a non-degenerate $ d \times d $ integer matrix
 whose eigenvalues are bigger than 1 in module,
$M^*$ is the conjugate  matrix to $M$,
$I_d$ denotes the unit $ d \times d$ matrix.
We  say that  numbers $k, n\in \zd $ are congruent modulo~ $M $
(write $k\equiv ~ n ~ \pmod {M} $) if $k-n=M\ell $, $ \ell\in\zd $.
The integer lattice $ \zd $ is split into cosets with
respect to the introduced relation of congruence.
The number of cosets is equal to
 $ | \det M | $ (see, e.g., \cite[\S 2.2]{NPS}).
Let us take  an arbitrary representative from each coset,
call them digits and denote the set of digits by $D (M) $.
Throughout the paper we  consider that such a matrix
$M$ is fixed, $m=|\det M|$,
$D (M) =\{s_0\ddd s_{m-1}\}$, $s_0=\nul$,
 $R (M) =\{M^{-1}s_0\ddd M^{-1}s_{m-1}\}$.

Let $\psi^{(\nu)}, \widetilde\psi^{(\nu)}\in L_2(\rd)$, $\nu-1\ddd r$.
Dual wavelet systems  $\{\psi_{jk}^{(\nu)}\}$, $\{\widetilde\psi_{jk}^{(\nu)}\}$, where
$\psi_{jk}^{(\nu)}:=m^{j/2}\psi^{(\nu)}(M^j\cdot+k)$,
 $\widetilde\psi_{jk}^{(\nu)}:=m^{j/2}\widetilde\psi^{(\nu)}(M^j\cdot+k)$,
$j,k\in\zd$, $\nu=1\ddd r$, form  dual frames if for each $f\in  L_2(\rd)$
\be
f=\sml
j{-\infty}\infty\sum\limits_{k\in\zd}\sml\nu1r \langle
f,\widetilde\psi_{jk}^{(\nu)}\rangle \psi_{jk}^{(\nu)}
\label{000}
\ee
and
\ba
A\|f\|^2<\sml j{-\infty}\infty\sum\limits_{k\in\zd}\sml\nu1r
|\langle f,\psi_{jk}^{(\nu)}\rangle|^2< B\|f\|^2,\ \ \  A,B>0,
\label{001}
\\
\widetilde A\|f\|^2<\sml j{-\infty}\infty\sum\limits_{k\in\zd}\sml\nu1r
|\langle f,\widetilde\psi_{jk}^{(\nu)}\rangle|^2<\widetilde B\|f\|^2, \ \ \ \widetilde A,\widetilde B>0.
\label{002}
\ea

Let a MRA in $L_2(\rd)$ be generated by a scaling function
$\phi$ which satisfies the refinement equation
$$
\widehat\phi(x)=m_0({M^*}^{-1}x) \widehat\phi({M^*}^{-1}x),
$$
where $m_0\in L_2(\Bbb T^d)$ is its mask (refinable mask).
For any  $m_\nu\in L_2(\Bbb T^d)$, there exists a unique
set of  functions $\mu_{\nu k}\in L_2(\Bbb T^d)$, $k=0\ddd r$,
(polyphase representatives of $m_\nu$) so that
\be
m_\nu(x)=\frac1{\sqrt m}\sml k0{m-1}\ex{(s_k,x)}\mu_{\nu k}(M^*x).
\label{2}
\ee
The functions $\mu_{\nu k}$ can be expressed by
$$
\mu_{\nu k}(x)=\frac1{\sqrt m}\sum\limits_{s\in D(M^*)}\exm{(M^{-1}s_k,x+s)}m_\nu({M^*}^{-1}(x+s)).
$$
It is clear from these formulas that a function $m_\nu$ is differentiable ($n$ times)  on $R(M^*)$
if and only if its polyphase representatives $\mu_{\nu k}$, $k=0\ddd m-1$, are  differentiable ($n$ times)
at the origin and $m_\nu$ is a trigonometric polynomial if and only if its polyphase representatives,
are  a trigonometric polynomials.

Let now another MRA be generated by a scaling function
$\widetilde\phi$ with a mask  $\widetilde m_0$.  According to{\em  Unitary Extension
Principle}~\cite{RS}, to construct dual
wavelet frames one finds wavelet masks $m_\nu, \widetilde m_\nu$,
$\nu=1\ddd r$, $r\ge m-1$, so that the polyphase matrices
\ban
{\cal M}:=\left(%
\begin{array}{ccc}
  \mu_{00} & \dots& \mu_{0, m-1} \\
  \vdots & \ddots & \vdots\\
  \mu_{r,0} & \dots & \mu_{r,m-1}\\
\end{array}%
\right),\ \ \
\widetilde{\cal M}:=\left(%
\begin{array}{ccc}
  \widetilde\mu_{00} & \dots& \widetilde\mu_{0, m-1} \\
  \vdots & \ddots & \vdots\\
  \widetilde\mu_{r,0} & \dots & \widetilde\mu_{r,m-1}\\
\end{array}%
\right),
\ean
satisfy
\be
{\cal M^T}\overline{\widetilde{\cal M}}=I_m,
\label{61}
\ee
and define wavelet functions by
\ban
\widehat\psi^{(\nu)}(x)&=&m_\nu({M^*}^{-1}x) \widehat\phi({M^*}^{-1}x),
\\
\widehat{\widetilde\psi^{(\nu)}}(x)&=&\widetilde m_\nu({M^*}^{-1}x)
 \widehat{\widetilde\phi}({M^*}^{-1}x).
\ean
The corresponding dual wavelet systems are $\{\psi_{jk}^{(\nu)}\}$, $\{\widetilde\psi_{jk}^{(\nu)}\}$
is said to be generated by $\phi$, $\widetilde\phi$ (or MRA-based).

It is known that if ${\cal M}=\widetilde{\cal M}$ then $\{\psi_{jk}^{(\nu)}\}$
is a tight frame in $L_2(\rd)$.  If ${\cal M}, \widetilde{\cal M}$ are arbitrary matrixes satisfying~(\ref{61}),
under some additional assumptions on $\phi,  \widetilde\phi,
m_\nu, \widetilde m_\nu$ (see \cite{98}, \cite{b16}, \cite[\S 2.7]{NPS}),
we can state that $\{\psi_{jk}^{(\nu)}\}$, $\{\widetilde\psi_{jk}^{(\nu)}\}$ are  dual frames  in $L_2(\rd)$.

\begin{defi}
We say that a  wavelet system  $\{\psi_{jk}^{(\nu)}\}$
has vanishing moments up to order $\alpha$, $\alpha\in\zd_+$,
(
has {\em $VM_\alpha$ property} in the sequel),
if $D^\beta\widehat\psi^{(\nu)}(\nul)=0$, $\nu=1\ddd r$,
for all $\beta\in\zd_+$, $\beta\le\alpha$.
\end{defi}

Usually it is more useful to control univariate
order of vanishing moment property (for example, to apply  Taylor formula).

\begin{defi}
We say that a  wavelet system  $\{\psi_{jk}^{(\nu)}\}$
has vanishing moments up to order $n$, $n\in\z_+$,
(has {\em  $VM^n$ property} in the sequel)
if $D^\beta\widehat\psi^{(\nu)}(\nul)=0$, $\nu=1\ddd r$,
for all $\beta\in\zd_+$, $[\beta]\le n$.
\end{defi}

\subsection {2. Why vanishing moments are needed for frames?}

\begin{theo} Any  wavelet frame $\{\psi_{jk}^{(\nu)}\}$  has $VM^0$ property
provided all the functions $\widehat\psi^{(\nu)}$ are bounded and
continuous  at the origin.
\label{t002}
\end{theo}

{\tt Proof.} Let $\{\psi_{jk}^{(\nu)}\}$ be a wavelet frame in $L_2(\rd)$.
Assume that $\widehat \psi^{(\nu)}(\nul)\ne0$ for at least one $\nu$. Take a function
$f\in L_2(\rd)$ whose Fourier transform is  compactly supported.
Using standard technique, we have
$$
\sum\limits_{k\in\zd}
|\langle f,\psi_{jk}^{(\nu)}\rangle|^2=\int\limits_{\rd}|
\widehat\psi^{(\nu)}(M^{-j}x|^2|\widehat f(x)|^2\,dx+R_j,
$$
where
$$
R_j\le\int\limits_{\rd}\sum\limits_{l\in\zd, l\ne0}|\widehat f(x)|\,|\widehat f(x+M^jl)|
\,|\widehat\psi^{(\nu)}(M^{-j}x|\, |\widehat\psi^{(\nu)}(M^{-j}x+l)|\,dx.
$$
If $j$ is large enough positive integer, then the supports
of the functions $\widehat f$ and $\widehat f(\cdot+M^jl)$ are disjoint for all $l\in\zd$, $l\ne\nul$.
Hence
$$
\sum\limits_{k\in\zd}
|\langle f,\psi_{jk}^{(\nu)}\rangle|^2=\int\limits_{\rd}|\widehat\psi^{(\nu)}(M^{-j}x|^2|\widehat f(x)|^2\,dx
\too\limits_{j\to+\infty}|\widehat\psi^{(\nu)}(\nul)|^2\|f\|^2.
$$
This contradict to~(\ref{001}).$\Diamond$

Let $\{\psi_{jk}^{(\nu)}\}$,  $\{\widetilde\psi_{jk}^{(\nu)}\}$
be dual frames. Decompositions~(\ref{000}) is said to have approximation
order $n$ if
there exist $C>0$ and $\lambda>1$ such that
 for any function $f$ in the Sobolev space $W_2^n$
\be
\left\|f-\sum\limits_{i<j}\sum\limits_{k\in\zd}\sml\nu1r
\langle f,\widetilde\psi_{ik}^{(\nu)}\rangle \psi_{ik}^{(\nu)}\right\|_2\le
C\frac{\|f\|_{W_2^n}}{\lambda^{jn}}.
\label{95}
\ee
Since all spectrum of the operator $M^{-1}$ is located in the circle
 $ |x| \le r (M^{-1}) $,
where $r (M^{-1}): = \lim_{i \to \infty} {\| M^{-i} \|^{1/i}} $ is the
 spectral radius of $M^{-1} $, and there exists at least
 one point of the spectrum  on the boundary of the circle, (\ref{95})
 follows from
$$
\left\|f-\sum\limits_{i<j}\sum\limits_{k\in\zd}\sml\nu1r
\langle f,\widetilde\psi_{ik}^{(\nu)}\rangle \psi_{ik}^{(\nu)}\right\|_2\le
C\|f\|_{W_2^n}\sml ij\infty \|M^{-i}\|^n
$$
with $\lambda$ which is less than the minimal module of an eigenvalue of $M$.
For a diagonal matrix $M=cI_d$, (\ref{95}) holds with  $\lambda=c$.

\begin{theo}
Let $\{\psi_{jk}^{(\nu)}\}$,  $\{\widetilde\psi_{jk}^{(\nu)}\}$
be dual wavelet frames,
$$
|\psi^{(\nu)}(x)|, |\widetilde\psi^{(\nu)}(x)|\le \frac{C_1}{(1+|x|)^\gamma},\ \ \gamma>n+d,
$$
for all $\nu=1\ddd r$, and almost all $x\in\rd$.  If $\{\widetilde\psi_{jk}^{(\nu)}\}$ has
$VM^{n-1}$ property, then decompositions~(\ref{000}) have approximation order $n$.
\label{t001}
\end{theo}

We need the following auxiliary statement for the proof.

\begin {lem}
Let $\eta$ be a positive bounded function decreasing    on
$ [0, \infty) $ so that $\eta(|x|)$ is summable on $\rd$.
Then there exists a constant $K $ depending  on  $\eta$ and $d$ such that
\be
\sum\limits_{k\in\zd} \eta (|x+k |) \eta (|y+k |)
\le  K\eta\left(\frac {|x-y |} {8} \right)
\label {X}
\ee
for all $x, y\in\rd $.
\label {l01}
\end {lem}

{\tt Proof.}
First of all note that there exist constants $K_1$, $K_2$ such that
$$
\eta(|t|)\le K_1,\ \ \ \sum\limits_{k\in\zd}\eta (|t+k |)\le K_2 \ \ \forall t\in\rd.
$$
Since both the left and the right hand sides of~(\ref{X}) are invariant with respect
to the operation $(x,y)\to(x+l, y+l)$, $\l\in\zd$, we can assume that $|x|\le\frac{\sqrt d}{2}$.
If $|y|\le2\sqrt d$, then
$$
\sum\limits_{k\in\zd} \eta (|x+k |) \eta (|y+k |)\le K_1K_2\le
\frac{K_1K_2}{\eta(3\sqrt d)}\eta(|x-y|).
$$
Now let  $|y|\ge2\sqrt d$. Since $|x-y|\le2|y|$ and $|y|\ge4|x|$, we have
\ban
\sum\limits_{|k|\le\frac{|y|}{2}} \eta (|x+k |) \eta (|y+k |)\le
\sum\limits_{|k|\le\frac{|y|}{2}} \eta (|x+k |)\eta\lll\frac{|y|}{2}\rrr\le
K_2\eta\lll\frac{|y|}{2}\rrr\le
\\
K_2\eta\lll\frac{|x-y|}{4}\rrr,
\\
\sum\limits_{|k|\ge\frac{|y|}{2}} \eta (|x+k |) \eta (|y+k |)\le
\sum\limits_{|k|\ge\frac{|y|}{2}} \eta \lll\frac{|k|}{2}\rrr \eta (|y+k |)\le
K_2\eta\lll\frac{|y|}{4}\rrr\le
\\
K_2\eta\lll\frac{|x-y|}{8}\rrr.\Diamond
\ean

{\tt Proof of Theorem~\ref{t001}.} Let $f\in W_2^n$, $j\in\z_+$.
It follows from~(\ref{000}) that
\be
\left\|f-\sum\limits_{i<j}\sum\limits_{k\in\zd}\sml\nu1r
\langle f,\widetilde\psi_{ik}^{(\nu)}\rangle \psi_{ik}^{(\nu)}\right\|_2\le
\sum\limits_{i\ge j}\sml\nu1r\left\|\sum\limits_{k\in\zd}
\langle f,\widetilde\psi_{ik}^{(\nu)}\rangle \psi_{ik}^{(\nu)}\right\|_2.
\label{004}
\ee
Since $VM^n$ property is equivalent to
$$
\int\limits_{\rd}y^\alpha\widetilde\psi_{ik}^{(\nu)}(y)\,dy=0,\ \  \nu=1\ddd r, i\in\z, k\in\zd, \ \ \
\forall \alpha\in\zd, [\alpha]\le n,
$$
using Taylor formula, we have
\ban
\left|\langle f,\widetilde\psi_{ik}^{(\nu)}\rangle \right|=
\left|\int\limits_{\rd}f(y)\widetilde\psi_{ik}^{(\nu)}(y)\,dy\right|=
\hspace{3cm}
\\
\left|\int\limits_{\rd}\,dy\, \widetilde\psi_{ik}^{(\nu)}(y)
\lll\sml l0{n-1}\frac{1}{l!}\lll (y-x)_1\frac{\partial}{\partial x_1}
+\dots+(y-x)_d\frac{\partial}{\partial x_d}\rrr^lf(x)+\right.\right.
\\
\left.\left.\int\limits_0^1\frac{(1-t)^{n-1}}{(n-1)!}
\lll (y-x)_1\frac{\partial}{\partial x_1}
+\dots+(y-x)_d\frac{\partial}{\partial x_d}\rrr^{n}f(x+t(y-x))\,dt\rrr
\right|\le
\\
C_2\sum\limits_{\alpha\in\zd\atop[\alpha]= n}\int\limits_{\rd}\,dy\,|x-y|^n|\widetilde\psi_{ik}^{(\nu)}(y)|
\,\int\limits_0^1|D^\alpha f(x+t(y-x))|\,dt.
\ean
From this, due to Lemma~\ref{l01} and Cauchy-Bunyakovskii inequality, we obtain
\ban
\left\|\sum\limits_{k\in\zd}
\langle f,\widetilde\psi_{ik}^{(\nu)}\rangle \psi_{ik}^{(\nu)}\right\|_2^2\le
\hspace{5cm}
\\
C_2^2\int\limits_{\rd}\,dx\lll \int\limits_{\rd}\,dy
|x-y|^n\sum\limits_{k\in\zd}|\widetilde\psi_{ik}^{(\nu)}(y)\psi_{ik}^{(\nu)}(x)|
\,\int\limits_0^1\sum\limits_{\alpha\in\zd\atop[\alpha]= n}|D^\alpha f(x+t(y-x))|\,dt\rrr^2\le
\\
C_3\int\limits_{\rd}\,dx\lll\int\limits_{\rd}\,dy\int\limits_0^1\,dt
\frac{m^i|x-y|^n}{(1+|M^i(x-y)|)^\gamma}
\lll \sum\limits_{\alpha\in\zd\atop[\alpha]= n}|D^\alpha f(x+t(y-x))|^2\rrr^{\frac12}\rrr^2=
\\
C_3\int\limits_{\rd}\,dx\lll\int\limits_{\rd}\,du\int\limits_0^1\,dt
\frac{m^i|u|^n}{(1+|M^iu|)^\gamma}
\lll \sum\limits_{\alpha\in\zd\atop[\alpha]= n}|D^\alpha f(x+tu)|^2\rrr^{\frac12}\rrr^2\le
 \\
C_3\int\limits_{\rd}\,dx\lll\int\limits_{\rd}\,du\int\limits_0^1\,dt
\frac{m^i|u|^n}{(1+|M^iu|)^\gamma}
 \sum\limits_{\alpha\in\zd\atop[\alpha]= n}|D^\alpha f(x+tu)|^2\right.\cdot
 \hspace{2cm}
 \\
\left.\int\limits_{\rd}\,du\int\limits_0^1\,dt
\frac{m^i|u|^n}{(1+|M^iu|)^\gamma}
\rrr\le
C_3\lll \|f\|_{W_p^2}\int\limits_{\rd}
\frac{m^i|u|^n\,du}{(1+|M^iu|)^\gamma}\rrr^2\le\hspace{2cm}
\\
C_3\lll \|f\|_{W_p^2}\|M^{-i}\|^n\int\limits_{\rd}
\frac{|v|^n\,dv}{(1+|v|)^\gamma}\rrr^2\le\hspace{3cm}
\\
C_3\lll \|f\|_{W_p^2}\|M^{-i}\|^n\int\limits_{\rd}
\frac{dv}{(1+|v|)^{\gamma-n}}\rrr^2=
C_4\|f\|^2_{W_p^2}\|M^{-i}\|^{2n}.\hspace{1cm}
\ean
It remains to combine this estimation with~(\ref{004}).
$\Diamond$

The scheme of the proof of Theorem~\ref{t001} can be applied
for some other  approximation problems. In a similar way,
it is possible to estimate coefficients  of decomposition~(\ref{000}),
to find the order of approximation at a point for functions with some
special local properties, e.g., for the class (introduced by
Calder\'on and Zygmund~\cite{9}) of functions $f\in L(\rd)$ such that,
$$
\int\limits_{|x-x_0|\le h}|f(x)-P(x-x_0)|\,dx= o(h^\alpha),
$$
where $P$ is an algebraic polynomial.  Detailed consideration of these problems is out
of our today's interest. Our main goal is investigation of $VM^n$ property for frames
and development of methods for  their construction.

\subsection {3. Polyphase characterization of vanishing moments property}

Now we will consider only
wavelet systems  $\{\psi_{jk}^{(\nu)}\}$,
 $\{\widetilde\psi_{jk}^{(\nu)}\}$ which are generated by scaling functions $\phi,  \widetilde\phi$
whose masks $m_0,  \widetilde m_0$ are continuous at the origin and
$m_0(\nul)= \widetilde m_0(\nul)=1$.

Assume  that the functions $\widehat\phi,
m_1\ddd m_r$ have  derivatives up to order $\alpha$
at the origin. It easily follows from Leibniz formula that
$VM_\alpha$ property holds if and only if
\be
D^\beta(m_\nu({M^*}^{-1}x)\Big|_{x=\nul}=0, \nu=1\ddd r,\ \ \
\forall \beta\in\zd_+, \beta\le\alpha.
\label{36}
\ee
In the case $r=m-1$, there exist different criterions for vanishing moment.
It is known~\cite{1} how to describe
vanishing moment property in terms of linear identities for Fourier coefficients of
the dual refinable mask (so-called {\em sum rule}).  Some other descriptions of masks
providing $VM_\alpha$ property are  found in terms of zero-conditions~\cite{1}
and in terms of containment in a quotient ideal~\cite{11}. The following polyphase criterion
was given in~\cite{36}: $VM_\alpha$ property is valid for   $\{\psi_{jk}^{(\nu)}\}$
if and only if there exist complex numbers $\lambda_\gamma$,
$\gamma\in\zd_+$,  $\gamma\le\alpha$, such that $\lambda_\nul=1$,
\be
D^\beta\widetilde\mu_{0k}(\nul)=\frac1{\sqrt m}
\sum\limits_{\nul\le\gamma\le\beta}\lambda_\gamma
\lll\beta\atop\gamma\rrr(-2\pi i M^{-1}s_k)^{\beta-\gamma} \ \ \ \forall    \beta\in\zd_+,  \beta\le\alpha,
\label{0}
\ee
for each $k=0\ddd m-1$. The set of parameters $\lambda_\gamma$  in~(\ref0)
is unique, and  $\lambda_\gamma$ does not depend on $\alpha$ due to the following statement.
\begin{prop}\cite{36}
If~(\ref{0}) is valid for the polyphase representatives of $\widetilde m_0$, then
\be
\lambda_\beta=D^\beta\lll\widetilde m_0({M^*}^{-1}x)\rrr\Big|_{x=\nul}
\label{35}
\ee
for all   $\beta\in\zd_+$,  $\beta\le\alpha$.
\label{pr1}
\end{prop}

So, in the case $r=m-1$,  $VM_\alpha$ property for  $\{\psi_{jk}^{(\nu)}\}$ depends only on
$\widetilde m_0$, i.e. only the first raw of the matrix $\widetilde{\cal M}$ is responsible for
vanishing moments of wavelets generated by the matrix ${\cal M}$. In the case $r>m-1$,
$VM_\alpha$ property for  $\{\psi_{jk}^{(\nu)}\}$ depends
also on the way of construction of matrixes ${\cal M}$, $\widetilde{\cal M}$.
This may be illustrated by the following example.

Let $d=1$, $M=m=2$, $\mu_{00}=\mu_{01}=\widetilde\mu_{00}=\widetilde\mu_{01}\equiv\frac1{\sqrt2}$,
$$
{\cal M}=\widetilde{\cal M}=\left(%
\begin{array}{cc}
  \frac1{\sqrt2} &  \frac1{\sqrt2} \\
   \frac12 &  -\frac12 \\
   \frac12 &  -\frac12 \\
\end{array}%
\right),
\ \ \ \ \
{\cal M}^\prime=\left(%
\begin{array}{cc}
  \frac1{\sqrt2} &  \frac1{\sqrt2} \\
   \frac1{\sqrt2} & 0 \\
   0 &  \frac1{\sqrt2} \\
\end{array}%
\right),
\ \ \
\widetilde{\cal M}^\prime=\left(%
\begin{array}{cc}
  \frac1{\sqrt2} &  \frac1{\sqrt2} \\
  \frac1{\sqrt2} &  -\frac1{\sqrt2} \\
   -\frac1{\sqrt2} &  \frac1{\sqrt2} \\
\end{array}%
\right).
$$
Either of pairs ${\cal M},\widetilde{\cal M}$ and ${\cal M}^\prime,\widetilde{\cal M}^\prime$
satisfies~(\ref{61}).  The matrixes
${\cal M},\widetilde{\cal M}$ generate wavelet masks
$m_1(x)=m_2(x)=\widetilde m_1(x)=\widetilde m_2(x)=\frac1{2\sqrt2}-\frac1{2\sqrt2}\ex x$.
It is clear that $m_1(0)=m_2(0)=\widetilde m_1(0)=\widetilde m_2(0)=0$, i.e.
for the corresponding wavelet systems  $VM_0$  property is valid. The matrixes
${\cal M}^\prime,\widetilde{\cal M}^\prime$ generate wavelet masks
$ m^\prime_1(x)=\frac1{2},\ \   m^\prime_2(x)=\frac1{2}\ex x,\ \
\widetilde m^\prime_1(x)=\frac1{2}-\frac1{2}\ex x, \widetilde m^\prime_2(x)=-\frac1{2}+\frac1{2}\ex x$,
and we have $m^\prime_1(0)\ne0, m^\prime_2(0)\ne0$.

\begin{theo}
Let  $\alpha\in\zd_+$, $r\ge m-1$,
the functions $\mu_{\nu,k},\widetilde\mu_{\nu k}\in L_2(\Bbb T^d)$, $\nu, k=0\ddd r$,
 have  derivatives up to order $\alpha$
at the origin, the matrixes ${\cal N}:=\{\mu_{\nu k}\}_{\nu,k=0}^r$ and
$\widetilde{\cal N}:=\{\widetilde\mu_{\nu k}\}_{\nu,k=0}^r$  satisfy
\be
{\cal N}\overline{{\widetilde{\cal N}}^T}=I_{r+1};
\label1
\ee
and let masks $\widetilde m_0, m_1\ddd m_{m-1}$ be defined by~(\ref{2}).
Then condition~(\ref{36}) is valid if and only if
\item (a) there exist $\lambda_\gamma\in\Bbb C$,
$\gamma\in\zd_+$,  $\gamma\le\alpha$, such that $\lambda_\nul=1$ and (\ref{0})~holds for $k=0\ddd m-1$;
\item (b)
$D^\gamma\widetilde\mu_{0k}(\nul)=0$, $k=m\ddd r$
for all    $\gamma\in\zd_+$,  $\gamma\le\alpha$.
\label{t1}
\end{theo}

{\tt Proof.}
Suppose that~(\ref{36}) is valid.
We will prove $(a)$ and $(b)$ by induction on $\alpha$.
Check the initial step  for $\alpha=0$. Let $m_\nu(\nul)=0$,
$\nu=1\ddd r$. It follows form~(\ref2) that
\be
\sml k0{m-1}\mu_{\nu k}(\nul)=0,\ \ \ \nu=1\ddd r.
\label4
\ee
On the other hand, by~(\ref1),
$$
\sml k0r\overline{\widetilde\mu_{0k}(\nul)}\mu_{\nu k}(\nul)=0,\ \ \ \nu=1\ddd m-1.
$$
Because of linear independence of the vectors $(\mu_{\nu0}(\nul)\ddd\mu_{\nu,r}(\nul))\in\r^{r+1}$,
$\nu=1\ddd r$, there exists $\lambda$ so that
$$
\widetilde\mu_{00}(\nul)=\dots=\widetilde\mu_{0, m-1}(\nul)=\lambda,\ \
\widetilde\mu_{0m}(\nul)=\dots=\widetilde\mu_{0, r}(\nul)=0.
$$
Taking into account the condition $\widetilde m_0(\nul)=1$ which is equivalent to
$$
\frac1{\sqrt m}\lll\widetilde\mu_{\nu0}(\nul)+\dots+\widetilde\mu_{\nu,m-1}(\nul)\rrr=1,
$$
we obtain $\lambda=\frac1{\sqrt m}$.

For the inductive step we assume that~(\ref{36})  is valid for $\alpha>\nul$ and $(a), (b)$  holds
for all $\alpha^\prime\in\zd_+$, $\alpha^\prime<\alpha$.
So, due to Proposition~\ref{pr1},
there exist constants $\lambda_\gamma\in \Bbb C$, $\gamma\in\zd_+$, $\gamma<\alpha$
such that~(\ref0) holds for all $\beta<\alpha$. If $\gamma\in\zd_+$,
$\gamma<\alpha$, due to~(\ref2) and Leibniz formula, we have
\ba
\frac1{\sqrt m}\sum\limits_{\nul\le\beta\le\alpha-\gamma}
\lll\alpha-\gamma\atop\beta\rrr \sml k0{m-1}(2\pi iM^{-1}s_k)^{\alpha-\beta-\gamma}D^\beta \mu_{\nu k}(\nul)=
\nonumber
\\
D^{\alpha-\gamma}m_\nu({M^*}^{-1}x)\Big|_{x=\nul}=0.
\label5
\ea
It follows from~(\ref1) that
$$
\sml k0r\overline{\widetilde\mu_{0k}}\mu_{\nu k}= 0,\ \ \nu=1\ddd m-1.
$$
Differentiating this equality $\alpha$ times gives
$$
\sum\limits_{\nul\le\beta\le\alpha}
\lll\alpha\atop\beta\rrr \sml k0r\overline{D^{\alpha-\beta }\widetilde\mu_{0 k}(\nul)}
D^\beta \mu_{\nu k}(\nul)=0.
$$
Taking into account the inductive hypotheses, we have
\be
\sum\limits_{\nul\le\beta\le\alpha}
\lll\alpha\atop\beta\rrr \sml k0{m-1}\overline{D^{\alpha-\beta }\widetilde\mu_{0 k}(\nul)}
D^\beta \mu_{\nu k}(\nul)+\sml kmr\overline{D^{\alpha}\widetilde\mu_{0 k}(\nul)}
\mu_{\nu k}(\nul)=0.
\label3
\ee
Multiply~(\ref5) by $\lll\alpha\atop\alpha-\gamma\rrr\overline{\lambda_\gamma}$
and subtract from~(\ref3). After the same manipulation with each $\gamma\in\zd_+$,
$\gamma<\alpha$, we obtain
\ban
0=\sml kmr\overline{D^{\alpha}\widetilde\mu_{0 k}(\nul)}
\mu_{\nu k}(\nul)+\sum\limits_{\nul\le\beta\le\alpha}\lll\alpha\atop\beta\rrr \sml k0{m-1}
\overline{D^{\alpha-\beta }\widetilde\mu_{0 k}(\nul)}D^\beta \mu_{\nu k}(\nul)-
\\
\frac1{\sqrt m}\sum\limits_{\nul\le\gamma<\alpha}\lll\alpha\atop\alpha-\gamma\rrr\overline{\lambda_\gamma}
\sum\limits_{\nul\le\beta\le\alpha-\gamma}
\lll\alpha-\gamma\atop\beta\rrr \sml k0{m-1}(2\pi iM^{-1}s_k)^{\alpha-\beta-\gamma}D^\beta \mu_{\nu k}(\nul)=
\\
\sml kmr\overline{D^{\alpha}\widetilde\mu_{0 k}(\nul)}
\mu_{\nu k}(\nul)+
\sum\limits_{\nul\le\beta\le\alpha}\lll\alpha\atop\beta\rrr \sml k0{m-1}\lll
\overline{D^{\alpha-\beta }\widetilde\mu_{0 k}(\nul)}-\right.
\\
\left.\frac1{\sqrt m}\sum\limits_{\gamma\ne\alpha\atop\nul\le\gamma<\alpha-\beta}
\overline{{\lll\alpha-\gamma\atop\beta\rrr \lll\alpha\atop\alpha-\gamma\rrr}
{\lll\alpha\atop\beta\rrr^{-1} }\lambda_\gamma(-2\pi iM^{-1}s_k)^{\alpha-\beta-\gamma}}\right)
D^\beta \mu_{\nu k}(\nul).
\ean
From this, taking into account that
\be
{\lll\alpha-\gamma\atop\beta\rrr\lll\alpha\atop\alpha-\gamma\rrr}
{\lll\alpha\atop\beta\rrr ^{-1}}=\frac{(\alpha-\beta)!}{\gamma!(\alpha-\beta-\gamma)!}=
\lll\alpha-\beta\atop\gamma\rrr,
\label6
\ee
and using the inductive hypotheses,  the sum over $\beta$ is deduced to
\ban
\sum\limits_{\nul\le\beta\le\alpha}\lll\alpha\atop\beta\rrr \sml k0{m-1}\lll
\overline{D^{\alpha-\beta }\widetilde\mu_{0 k}(\nul)}-\right.\hspace{6cm}
\\
\left.\frac1{\sqrt m}\sum\limits_{\gamma\ne\alpha\atop\nul\le\gamma\le\alpha-\beta}
\overline{{\lll\alpha-\beta\atop\gamma\rrr \ }
\lambda_\gamma(-2\pi iM^{-1}s_k)^{\alpha-\beta-\gamma}}\right)
D^\beta \mu_{\nu k}(\nul)=
\\
\sml k0{m-1}\lll
\overline{D^{\alpha}\widetilde\mu_{0 k}(\nul)-
\frac1{\sqrt m}\sum\limits_{\nul\le\gamma<\alpha}
{\lll\alpha\atop\gamma\rrr \ }
\lambda_\gamma(-2\pi iM^{-1}s_k)^{\alpha-\gamma}}\right)
\mu_{\nu k}(\nul).
\ean
So, we have
\ban
\sml k0{m-1}\lll
\overline{D^{\alpha}\widetilde\mu_{0 k}(\nul)-
\frac1{\sqrt m}\sum\limits_{\nul\le\gamma<\alpha}
{\lll\alpha\atop\gamma\rrr \ }
\lambda_\gamma(-2\pi iM^{-1}s_k)^{\alpha-\gamma}}\right)
\mu_{\nu k}(\nul)+
\\
\sml kmr\overline{D^{\alpha}\widetilde\mu_{0 k}(\nul)}
\mu_{\nu k}(\nul)=0.
\ean
Similarly to the arguments for the initial step, it follows from~(\ref4) that
there exists $\lambda_\alpha$ such that
\ban
D^{\alpha }\widetilde\mu_{0 k}(\nul)-
\frac1{\sqrt m}\sum\limits_{\nul\le\gamma<\alpha}
{\lll\alpha\atop\gamma\rrr
\lambda_\gamma(-2\pi iM^{-1}s_k)^{\alpha-\gamma}}=\frac{\lambda_\alpha}{\sqrt m}, \ \ \ k=0\ddd m-1,
\\
D^{\alpha }\widetilde\mu_{0 k}(\nul)=0, \ \ \ k=m\ddd r.
\ean
Thus, (\ref0) is valid for $\beta=\alpha$ as was to be proved.

Now we assume that $(a), (b)$ are valid.
We  will  prove~(\ref{36})  by induction on~$\alpha$.
If~(\ref0) is valid for $\alpha=\nul$, then $\widetilde\mu_{0 k}(\nul)=1/\sqrt m$,
$k=0\ddd m-1$. It follows from~(\ref{1}) and $(b)$ that
$$
\mu_{\nu 0}(\nul)+\dots+\mu_{\nu, m-1}(\nul)=0,\ \ \nu=1\ddd r.
$$
Hence, on the basis of~(\ref{2}), $m_\nu(\nul)=0$, $\nu=1\ddd r$, what proves
the initial step.

For the inductive step, we assume that $(a), (b)$ is valid for $\alpha>\nul$ and~(\ref{36})  holds
for all $\alpha^\prime\in\zd_+$, $\alpha^\prime<\alpha$, i.e.
$$
D^{\alpha-\gamma}m_\nu({M^*}^{-1}x)\Big|_{x=\nul}=0,\ \
\gamma\in\zd_+, \ \ \gamma\ne\nul,  \ \ \gamma\le\alpha.
$$
This yields~(\ref5) for $\gamma\ne\nul$. Multiply~(\ref5)
by $\lll\alpha\atop\alpha-\gamma\rrr\overline{\lambda_\gamma}$
and add to~(\ref2) differentiated $\alpha$ times. After the same manipulation
with each $\gamma\in\zd_+$, $\gamma<\alpha$, we obtain
\ban
D^{\alpha}m_\nu({M^*}^{-1}x)\Big|_{x=\nul}=
\frac1{\sqrt m}\sum\limits_{\nul\le\beta\le\alpha}
\lll\alpha\atop\beta\rrr \sml k0{m-1}(2\pi iM^{-1}s_k)^{\alpha-\beta}D^\beta \mu_{\nu k}(\nul)+
\\
\frac1{\sqrt m}\sum\limits_{\nul<\gamma\le\alpha}\lll\alpha\atop\alpha-\gamma\rrr\overline{\lambda_\gamma}
\sum\limits_{\nul\le\beta\le\alpha-\gamma}
\lll\alpha-\gamma\atop\beta\rrr \sml k0{m-1}(2\pi iM^{-1}s_k)^{\alpha-\beta-\gamma}D^\beta \mu_{\nu k}(\nul)=
\\
\frac1{\sqrt m}\sum\limits_{\nul\le\beta\le\alpha}\lll\alpha\atop\beta\rrr
\sum\limits_{\nul\le\gamma\le\alpha-\beta}\overline{\lambda_\gamma}\lll\alpha\atop\alpha-\gamma\rrr
\lll\alpha-\gamma\atop\beta\rrr\lll\alpha\atop\beta\rrr^{-1}\cdot
\\
\sml k0{m-1}(2\pi iM^{-1}s_k)^{\alpha-\beta-\gamma}D^\beta \mu_{\nu k}(\nul).
\ean
Due to~(\ref6), (\ref0) and $(b)$, this yields
\ban
D^{\alpha}m_\nu({M^*}^{-1}x)\Big|_{x=\nul}=\hspace{8cm}
\\
\frac1{\sqrt m}\sum\limits_{\nul\le\beta\le\alpha}\lll\alpha\atop\beta\rrr\sml k0{m-1}
\overline{\sum\limits_{\nul\le\gamma\le\alpha-\beta}\lambda_\gamma\lll\alpha-\beta\atop\gamma\rrr
(-2\pi iM^{-1}s_k)^{\alpha-\beta-\gamma}}D^\beta \mu_{\nu k}(\nul)=
\\
\sum\limits_{\nul\le\beta\le\alpha}\lll\alpha\atop\beta\rrr
\sml k0{m-1}\overline{D^{\alpha-\beta}\widetilde\mu_{0k}(\nul)}D^\beta \mu_{\nu k}(\nul)=
D^\alpha\lll\sml k0{m-1}\overline{\widetilde\mu_{0k}(x)} \mu_{\nu k}(x)\rrr\Bigg|_{x=\nul}=
\\
D^\alpha\lll\sml k0r\overline{\widetilde\mu_{0k}(x)} \mu_{\nu k}(x)\rrr\Bigg|_{x=\nul}.
\ean
It follows from~(\ref1) that
$
D^{\alpha}m_\nu({M^*}^{-1}x)\Big|_{x=\nul}=0
$
as was to be proved. $\Diamond$

It is not difficult to see that $\{\psi_{jk}^{(\nu)}\}$ has  $VM^n$ property if and only if
(\ref{36}) is valid for all $\alpha\in\zd$, $[\alpha]\le n$.
The following statement follows immediately from Theorem~\ref{t1} and Proposition~\ref{pr1}.

\begin{theo} Let  $n\in\z_+$, \ \ $r\ge m-1$, the functions $\mu_{\nu,k},\ \widetilde\mu_{\nu k}\in L_2(\Bbb T^d)$,
$\nu, k=0\ddd r$,
 have derivatives up to order $n$
at the origin, the matrixes \linebreak ${\cal N}:=\{\mu_{\nu k}\}_{\nu,k=0}^r$,
$\widetilde{\cal N}:=\{{\widetilde \mu_{\nu k}}\}_{\nu,k=0}^r$ satisfy~(\ref1);
and let masks $\widetilde m_0, m_1\ddd m_{m-1}$ be defined by~(\ref{2}).
Then condition~(\ref{36}) is valid for all $\alpha\in\zd$, $[\alpha]\le n$, if and only if
\item (a) there exist  $\lambda_\gamma\in\Bbb C$,
$\gamma\in\zd_+$,  $[\gamma]\le n$, such that $\lambda_\nul=1$ and (\ref{0})~holds for $k=0\ddd m-1$;
\item (b)
$D^\gamma\widetilde\mu_{0k}(\nul)=0$, $k=m\ddd r$
for all    $\gamma\in\zd_+$,  $[\gamma]\le n$.
\label{t2}
\end{theo}

Let  $n\in\z_+$, we will denote by $L_\infty^{(n)}$ the class of complex-valued functions
which are in $L_\infty(\Bbb T^d)$ and have continuous derivatives up to order $n$ at the origin.

\begin{lem}
Let $\mu_{\nu 0},\widetilde\mu_{\nu 0}\in L_\infty^{(n)}$, $\nu=0\ddd r$, and
\be
\sml \nu0r\mu_{\nu 0}\overline{\widetilde\mu_{\nu 0}}= 1.
\label{55}
\ee
Then there exist
functions $\mu_{\nu k}, \widetilde\mu_{\nu k}\in L_\infty^{(n)}$, $\nu=0\ddd r$, $k=1\ddd r$,
such that
\be
\sml \nu0r\mu_{\nu l}\overline{\widetilde\mu_{\nu k}}= \delta_{kl},\ \ \ k,l=0\ddd r.
\label{54}
\ee
\label{l1}
\end{lem}

{\tt Proof.} Set
$$
\mu^\prime_{\nu 0}:=
\left\{
\begin{array}{ll}
\mu_{\nu 0}\Bigg/\sqrt{\sml l0r|\mu_{l0}|^2},& \mbox{if}\ \ \  \sqrt{\sml l0r|\mu_{l0}|^2}\ne0,
\\
1/\sqrt m,& \mbox{if}\ \ \  \sqrt{\sml l0r|\mu_{l0}|^2}=0,
\end{array}
\right.
\ \ \ \nu=0\ddd r.
$$
It is clear that the functions $\mu^\prime_{\nu 0}$ are essentially bounded and
\be
\sml \nu0r|\mu^\prime_{\nu 0}|^2= 1.
\label{56}
\ee
It follows from~(\ref{55}) that $\sml l0r|\mu_{l0}(\nul)|^2\ne0$.
So, $\mu^\prime_{\nu 0}\in L_\infty^{(n)}$, $\nu=0\ddd r$.
Let us extend the unit vector $\mu^\prime_{0 0}\ddd \mu^\prime_{r 0}$ to a
unitary matrix. Due to~(\ref{56}), there exist $\nu_0$ so that $\mu^\prime_{\nu_0 0}(\nul)\ne1$.
We may assume that $\nu_0=0$ (else we will interchange $\mu^\prime_{\nu_0 0}$ and $\mu^\prime_{0 0}$,
extend this new vector to a unitary matrix and  interchange its $0$-th and $\nu_0$-th rows).
Due to Householder transform, an extension to a unitary matrix may be realized by:
$$
\mu^\prime_{0 k}=\overline{\mu^\prime_{k 0}}\frac{1-\mu^\prime_{00}}{1-\overline{\mu^\prime_{00}}},\ \
\mu^\prime_{\nu k}=\delta_{lk}-
\frac{\mu^\prime_{\nu 0}\overline{\mu^\prime_{k 0}}}{1-\overline{\mu^\prime_{00}}}, \ \ \ \nu, k=1\ddd r.
$$
Because of~(\ref{56}), we have $|\mu^\prime_{\nu 0}|\le\sqrt{1-|\mu^\prime_{00}|^2}$, $\nu=1\ddd r$.
This yields essential boundedness of the functions $\mu^\prime_{\nu k}$. Since
$1-\mu^\prime_{00}(\nul)\ne0$, it follows that $\mu^\prime_{\nu k}\in L_\infty^{(n)}$, $\nu,k=0\ddd r$.
Set
\ban
&&\widetilde\mu_{\nu k}:=\mu^\prime_{\nu k}, \ \ \ \nu=0\ddd r, \ k=1\ddd, r,
\\
&&\widetilde Q_k:=(\widetilde \mu_{0k}\ddd\widetilde \mu_{rk}), k=0\ddd r,
\\
&&Q_0:=(\mu_{00}\ddd\mu_{r0}),
Q_k:=\widetilde Q_k-\widetilde Q_k\overline{\widetilde Q^T_0}Q_0,\ \ \  k=1\ddd r.
\ean
It is not difficult to see that the entries of $Q_k$ are in $ L_\infty^{(n)}$ and
$Q_k\overline{\widetilde Q_l^T}=\delta_{kl}$, \ \ \ $k,l=0\ddd r$.
It remains to denote by $\mu_{\nu k}$ the $\nu$-th component of $Q_k$.$\Diamond$

\begin{lem}
Let $A$ be a class of complex-valued functions such that
\item  (i) if
$f,g\in A$, $a,b\in{\Bbb C}$ then $ af+bg\in A$,
\item  (ii) if  $f,g\in A$, then  $fg\in A$,\newline
and let $\cal A$ be a class of matrixes whose entries
are in $A$.
If any two $n\times 1$ matrixes $Q$, $\widetilde
Q\in{\cal A}$ satisfying $Q^T\overline{\widetilde{Q}}=1$  can be
extended to $n\times n$ matrixes ${\cal N}, \widetilde{\cal
N}\in{\cal A}$ satisfying ${\cal N}^T\overline{\widetilde{\cal
N}}=I_n$, then any  two $n\times j$  matrixes $\cal M$,
$\widetilde {\cal M}\in{\cal A}$, $1<j<n$,  satisfying
$ {\cal M}^T\overline{\widetilde{ {\cal M}}}=I_j$.
 can be extended to $n\times n$ matrixes ${\cal N}, \widetilde{\cal N}\in{\cal A}$ satisfying
${\cal N}^T\overline{\widetilde{\cal N}}=I_n$,
  \label{l2}
\end{lem}

{\tt Proof.} We will prove by induction on $j$. The base for $j=1$
is given. Let us check the inductive step $j-1\to j$. Let $j\times n$ matrixes
$\cal M, \widetilde {\cal M}\in \cal A$ satisfy $ {\cal M}^T\overline{\widetilde{ {\cal M}}}=I_j$. Denote by
$Q_k$, $\widetilde Q_k$ the $k$-th columns respectively  of
$\cal M$, $\widetilde {\cal M}$. Due to the statement of the
base, the matrixes $Q_1$, $\widetilde Q_1\in{\cal A}$  can be extended to
$n\times n$ matrixes ${{\cal N^\prime}}, \widetilde{\cal N^\prime}\in{\cal A}$
satisfying ${\cal N^\prime}^T\overline{\widetilde{\cal N^\prime}}=I_n$,
 Let $Q^\prime_k$, $\widetilde Q^\prime_k$, $k=2\ddd n$, denote
  the $k$-th columns respectively  of
  ${{\cal N^\prime}}, \widetilde{\cal N^\prime}\in{\cal  A}$.
Fix a point $x$ for which ${Q^\prime_l(x)}^T\overline{\widetilde Q^\prime}_k(x)=\delta_{kl}$.
Since the vectors $Q^\prime_2(x)\ddd Q^\prime_n(x)$
form a basis for the orthogonal complement to $\widetilde Q^\prime_1(x)$
in $\r^n$, we have
$$
Q_k(x)=\sml l2n\alpha_{lk}(x)Q^\prime_l(x),\ \ \ k=2\ddd j.
$$
Similarly,
$$
\widetilde Q_k(x)=\sml l2n\widetilde\alpha_{lk}(x)\widetilde Q^\prime_l(x),\ \ \ k=2\ddd j.
$$
It is clear that $\alpha_{lk}, \widetilde\alpha_{lk}\in A$ and
$\sml l2n \alpha_{lk}\overline{\widetilde\alpha}_{lk^\prime}=\delta_{k k^\prime}$, $k, k^\prime=2\ddd j$.
Due to the inductive hypotheses, there exist functions $\alpha_{lk}, \widetilde\alpha_{lk}\in A$,
$l=2\ddd n$, $k=j+1\ddd n$, such that
\be
\sml l2n \alpha_{lk}\overline{\widetilde\alpha}_{lk^\prime}=\delta_{k k^\prime},\ \ \ k, k^\prime=2\ddd n.
\label{53}
\ee
Set
$$
Q_k:=\sml l2n\alpha_{lk}Q^\prime_l,\ \ \
\widetilde Q_k:=\sml l2n\widetilde\alpha_{lk}\widetilde Q^\prime_l,\ \ \
 k=j+1\ddd n.
$$

Because of~(\ref{53}) and  biorthogonality of the systems $Q_1, {Q^\prime_2}\ddd {Q^\prime_n}$ and
${\widetilde Q_1},{\widetilde Q^\prime_2}\ddd {\widetilde Q^\prime_n}$,
we obtain
\ban
{Q_l(x)}^T\overline{\widetilde Q}_k(x)=\delta_{kl}, \ \ k,l=1\ddd n.
\ean
To complete the proof it remains to introduce matrixes ${\cal N}$ and  $\widetilde{\cal N}$
whose columns are respectively $Q_1\ddd Q_n$ and  $\widetilde Q_1\ddd \widetilde Q_n$. $\Diamond$

Now we are ready to give a necessary condition for $VM^n$ property.

\begin{theo}
Let dual wavelet systems $\{\psi_{jk}^{(\nu)}\}$,
$\{\widetilde\psi_{jk}^{(\nu)}\}$ be generated by refinable
functions $\phi,  \widetilde\phi$ whose Fourier transforms have
derivatives up to order $n$ at the origin, and let the entries
$\mu_{\nu k},\widetilde \mu_{\nu k}$ of the corresponding
 polyphase matrixes ${\cal M}, \widetilde{\cal M}$ be in $L^{(n)}_\infty$.
 If $\{\psi_{jk}^{(\nu)}\}$ has  $VM^n$ property, then
 there exist complex numbers $\lambda_\gamma$,
$\gamma\in\zd_+$,  $[\gamma]\le n$, $\lambda_\nul=1$ such
that~(\ref{0})  holds  for $k=0\ddd m-1$.
\label{t3}
\end{theo}

{\tt Proof.} Due to Lemmas~\ref{l1}, \ref{l2}, the $(r+1)\times m$
matrixes ${\cal M}, \widetilde{\cal M}$ can be extended to
$(r+1)\times (r+1)$ matrixes ${\cal N}, \widetilde{\cal N}$ such
that their entries are in $L^{(n)}_\infty$ and ${\cal
N}^T\overline{\widetilde{\cal N}}=I_{r+1}$. It remains to apply
conditions (a) of Theorem~\ref{t2}.$\Diamond$

Next let us discuss how to construct a dual wavelet system with
$VM^n$ property generated by a pair of refinable functions $\phi,
\widetilde\phi$ with  masks $m_0. \widetilde m_0$ respectively.
Due to Theorem~\ref{t3}, we know that the polyphase
representatives of $\widetilde m_0$ should satisfy~(\ref{0}). But
this theorem said nothing about $m_0$. Consider the following
example. Let $d=1$, $M=m=2$, $s_0=0, s_1=1$,
 $\mu_{00}=\widetilde\mu_{00}=\mu_{01}\equiv\frac1{\sqrt2},\ \ \
  \widetilde\mu_{01}=\frac1{\sqrt2}\lll1-\frac i2\sin2\pi x\rrr$. It is clear that~(\ref{0})
  holds for $n=1$ (with $\lambda_0=1$,
$\lambda_1=0$). Assume that there exists  dual wavelet systems
$\{\psi_{jk}^{(\nu)}\}$, $\{\widetilde\psi_{jk}^{(\nu)}\}$ with
$VM^1$ property (for $\{\psi_{jk}^{(\nu)}\}$) generated by these functions. This means that
the polyphase matrixes ${\cal M}, \widetilde{\cal M}$, whose first rows
are  respectively $(\mu_{00}, \mu_{01})$, $(\widetilde\mu_{00},
\widetilde\mu_{01})$, satisfy~(\ref{61}) and such
that~(\ref{36}) holds for the corresponding wavelet masks. Due to
Lemmas~\ref{l1}, \ref{l2}, the  matrixes ${\cal M},
\widetilde{\cal M}$ can be extended to  matrixes ${\cal N},
\widetilde{\cal N}$ such that their entries $\mu_{\nu k}$,
$\widetilde\mu_{\nu k}$, $\nu, k=0\ddd r$, are in $L^{(2)}_\infty$
and ${\cal N}^T, \overline{\widetilde{\cal N}}$ are mutually
inverse. It follows from Theorem~\ref{t2} that
\be
\frac{d}{dx}\lll\sml k2r\mu_{0 k}(x)\overline{\widetilde\mu_{0
k}(x)}\rrr\Bigg|_{x=0}=0. \label{62}
\ee
But
\ban
\frac{d}{dx}\lll\sml k2r\mu_{0 k}(x)\overline{\widetilde\mu_{0
k}(x)}\rrr\Bigg|_{x=0}= \frac{d}{dx}\lll1-\sml k01\mu_{0
k}(x)\overline{\widetilde\mu_{0 k}(x)}\rrr\Bigg|_{x=0}=
\\
\frac{d}{dx}\lll1-\frac12-\frac12\lll1-\frac i2\sin2\pi x\rrr\rrr\Bigg|_{x=0}=
\frac{d}{dx}\lll\frac i4\sin2\pi x\rrr\Bigg|_{x=0}\ne0.
\ean

So, we see that a generating refinable function $\phi$ should be
also chosen properly to provide  $VM^n$ property for
$\{\psi_{jk}^{(\nu)}\}$.

The following statement gives a sufficient condition for a pair of
refinable functions to generate  dual wavelet systems with  $VM^n$
property.

\begin{theo}
Let $\phi,  \widetilde\phi$ be refinable functions, their Fourier
transforms $\widehat\phi,  \widehat{\widetilde\phi}$  have
derivatives up to order $n$ at the origin, $\widehat\phi(\nul)=
\widehat{\widetilde \phi}(\nul)=1$, and let $\mu_{00}\ddd\mu_{0,
m-1}$, $\widetilde \mu_{00}\ddd\widetilde \mu_{0, m-1}\in
L^{(n)}_\infty$ be the  polyphase representatives of their masks.
If  there exist complex numbers $\lambda_\gamma$,
$\gamma\in\zd_+$,  $[\gamma]\le n$, $\lambda_\nul=1$, such that
(\ref{0})  holds  for $k=0\ddd m-1$  and there exist functions
 $\mu_{0k}$, $\widetilde \mu_{0k}\in L^{(n)}_\infty$, $k=m\ddd r$, such that
 \ban
 &&\sml k0r\mu_{0 k}\overline{\widetilde\mu_{0k}}=1,
 \\
&&D^\beta\widetilde\mu_{0k}(\nul)=0,\ k=m\ddd r,\ \ \ \forall
\beta\in\zd_+, [\beta]\le n.
 \ean
then the  functions  $\phi,  \widetilde\phi$ generate  dual
wavelet systems $\{\psi_{jk}^{(\nu)}\}$,
$\{\widetilde\psi_{jk}^{(\nu)}\}$ with  $VM^n$ property for
$\{\psi_{jk}^{(\nu)}\}$.
\label{t4}
\end{theo}

{\tt Proof.} Set $Q=(\mu_{00}\ddd\mu_{0r})$, $\widetilde
Q=(\widetilde\mu_{00}\ddd\widetilde\mu_{0r})$. Due to
Lemma~\ref{l1}, the $1\times (r+1)$ matrixes $Q, \widetilde Q$ can
be extended to $(r+1)\times (r+1)$ matrixes ${\cal N}=\{\mu_{\nu
k}\}_{\nu,k=0}^r, \widetilde{\cal N}=\{\widetilde\mu_{\nu
k}\}_{\nu,k=0}^r$ such that their entries are in $L^{(n)}_\infty$
and ${\cal N}\overline{{\widetilde{\cal N}}^T}=I_{r+1}$. So, the
matrixes \ban
{\cal M}:=\left(%
\begin{array}{ccc}
  \mu_{00} & \dots& \mu_{0, m-1} \\
  \vdots & \ddots & \vdots\\
  \mu_{r,0} & \dots & \mu_{r,m-1}\\
\end{array}%
\right),\ \ \
\widetilde{\cal M}:=\left(%
\begin{array}{ccc}
  \widetilde\mu_{00} & \dots& \widetilde\mu_{0, m-1} \\
  \vdots & \ddots & \vdots\\
  \widetilde\mu_{r,0} & \dots & \widetilde\mu_{r,m-1}\\
\end{array}%
\right), \ean satisfy~(\ref{61}). It follows from Theorem~\ref{t2}
that the corresponding wavelet masks $m_1\ddd m_{m-1}$
satisfy~(\ref{36}) for all $\alpha\in\zd$, $[\alpha]\le n$, what
was to be proved. $\Diamond$

Applied mathematicians and engineers are especially interested in
construction of compactly supported wavelet systems. To provide
this property  generating refinable functions should be  compactly
supported and wavelet masks should be trigonometric polynomials.

\begin{theo}
Let $\phi,  \widetilde\phi$ be  compactly supported refinable
functions with polynomial masks,
 $\widehat\phi(\nul)= \widehat{\widetilde \phi}(\nul)=1$,
and let $\mu_{00}\ddd\mu_{0, m-1}$, $\widetilde
\mu_{00}\ddd\widetilde \mu_{0, m-1}$ be the  polyphase
representatives of their masks. If  there exist complex numbers
$\lambda_\gamma$, $\gamma\in\zd_+$, $[\gamma]\le n$, such that~(\ref{0})
 holds  for $k=0\ddd m-1$  and there exist trigonometric polynomials
 $\mu_{0k}$, $\widetilde \mu_{0k}$, $k=m\ddd r$, such that
 \ban
 &&\sml k0r\mu_{0 k}\overline{\widetilde\mu_{0k}}=1,
 \\
&&D^\beta\widetilde\mu_{0k}(\nul)=0,\
k=m\ddd r,\ \ \ \forall \beta\in\zd_+, [\beta]\le n.
 \ean
then the  functions  $\phi,  \widetilde\phi$ generate dual
compactly supported  wavelet systems $\{\psi_{jk}^{(\nu)}\}$,
$\{\widetilde\psi_{jk}^{(\nu)}\}$ with  $VM^n$ property for
$\{\psi_{jk}^{(\nu)}\}$.
\label{t5}
\end{theo}

{\tt Proof.} Set $Q=(\mu_{00}\ddd\mu_{0r})$, $\widetilde
Q=(\widetilde\mu_{00}\ddd\widetilde\mu_{0r})$. Due to Suslin's
solution of a generalized  Serre conjecture\cite{3}, the row  $Q$
can be extended to a unimodular matrix with polynomial entries.
After this it is not difficult to find $(r+1)\times (r+1)$
matrices ${\cal N}, \widetilde{\cal N}$ extending
 $Q, \widetilde Q$. such that their entries are trigonometric polynomials and
${\cal N}\overline{{\widetilde{\cal N}}^T}=I_{r+1}$
(see~\cite{40},~\cite{32}, \cite[\S 2.6]{NPS}). Next we repeat the
arguments of the previous proof. $\Diamond$

Let us return to the example before Theorem~\ref{t4}. We could not succeed
with  $VM^1$ property because the derivative of $\sml k01\mu_{0
k}(x)\overline{\widetilde\mu_{0 k}(x)}$ did not vanish at the zero.
Try to change $\mu_{01}\equiv\frac1{\sqrt2}$ for $\mu_{01}=\widetilde\mu_{01}$.
Now we have
\ban
\frac{d}{dx}\lll\sml k01\mu_{0
k}(x)\overline{\widetilde\mu_{0 k}(x)}\rrr\Bigg|_{x=0}=
\frac{d}{dx}\lll\frac 18\sin^22\pi x\rrr\Bigg|_{x=0}=0.
\ean
So, providing  condition~(\ref0) for both the masks $m_0, \widetilde m_0$
improved the situation. Taking into account this observation, let us  consider
generating masks $m_0, \widetilde m_0$ whose polyphase
representatives satisfy
\ba
D^\beta\mu_{0k}(\nul)=\frac1{\sqrt m}
\sum\limits_{\nul\le\gamma\le\beta}\lambda_\gamma
\lll\beta\atop\gamma\rrr(-2\pi i M^{-1}s_k)^{\beta-\gamma} \ \ \ \forall    \beta\in\zd_+,  [\beta]\le n,
\label{58}
\\
D^\beta\widetilde\mu_{0k}(\nul)=\frac1{\sqrt m}
\sum\limits_{\nul\le\gamma\le\beta}\widetilde\lambda_\gamma
\lll\beta\atop\gamma\rrr(-2\pi i M^{-1}s_k)^{\beta-\gamma} \ \ \ \forall    \beta\in\zd_+,  [\beta]\le n.
\label{59}
\ea

\begin{theo}
Let dual wavelet systems $\{\psi_{jk}^{(\nu)}\}$, $\{\widetilde\psi_{jk}^{(\nu)}\}$
be generated by refinable functions $\phi,  \widetilde\phi$ whose  Fourier transforms have  derivatives up
to order $n$ at the origin,  the entries  $\mu_{\nu k},\widetilde \mu_{\nu k}$ of the corresponding
 polyphase matrixes ${\cal M}, \widetilde{\cal M}$ be in $L^{(n)}_\infty$ and
 there exist complex numbers $\lambda_\gamma, \widetilde\lambda_\gamma$,
$\gamma\in\zd_+$,  $[\gamma]\le n$, $\lambda_\nul=\widetilde\lambda_\nul=1$,
such that~(\ref{58}), (\ref{59}) are fulfilled  for $k=0\ddd m-1$. If at least one of the systems
$\{\psi_{jk}^{(\nu)}\}$, $\{\widetilde\psi_{jk}^{(\nu)}\}$ has  $VM^n$ property,
then
\ba
\sum\limits_{\nul\le\gamma\le \alpha}\lll\alpha\atop\gamma\rrr
\lambda_\gamma \overline{\widetilde\lambda_{\alpha-\gamma}}=0 \ \ \
\forall \alpha\in\zd_+, 0<[\alpha]\le n.
\label{57}
\ea
\label{t6}
\end{theo}

{\tt Proof.}
Let $\alpha\in\zd_+$, $0<[\alpha]\le n$, $k=0\ddd m-1$, $\rho:=2\pi M^{-1}s_k$.
Due to~(\ref{58}), (\ref{59}), we have
\ba
m D^\alpha\lll \mu_{0 k}(x)\overline{\widetilde\mu_{0 k}(x)} \rrr\Big|_{x=\nul}=
\sum\limits_{\nul\le\beta\le \alpha}\lll\alpha\atop\beta\rrr
D^\beta\mu_{0k}(\nul)\overline{D^{\alpha-\beta}\widetilde\mu_{0k}(\nul)}=\nonumber
\\
\sum\limits_{\nul\le\beta\le \alpha}\lll\alpha\atop\beta\rrr
\sum\limits_{\nul\le\gamma\le\beta}
\lll\beta\atop\gamma\rrr\lambda_{\beta-\gamma}(-i\rho)^\gamma
\sum\limits_{\nul\le\delta\le\alpha-\beta}
\lll\alpha-\beta\atop\delta\rrr\overline{\widetilde\lambda_{\alpha-\beta-\delta}}(i\rho)^\delta=\nonumber
\\
\sum\limits_{\nul\le\gamma\le\alpha}\sum\limits_{\gamma\le\beta\le \alpha}\sum\limits_{\nul\le\delta\le\alpha-\beta}
\lll\alpha\atop\beta\rrr\lll\beta\atop\gamma\rrr\lll\alpha-\beta\atop\delta\rrr
(-\odin)^\gamma(i\rho)^{\gamma+\delta}\lambda_{\beta-\gamma}
\overline{\widetilde\lambda_{\alpha-\beta-\delta}}=\nonumber
\\
\sum\limits_{\nul\le\gamma\le\alpha}\sum\limits_{\nul\le\delta\le\alpha-\gamma}
\sum\limits_{\gamma\le\beta\le \alpha-\delta}
\lll\alpha\atop\beta\rrr\lll\beta\atop\gamma\rrr\lll\alpha-\beta\atop\delta\rrr
(-\odin)^\gamma(i\rho)^{\gamma+\delta}\lambda_{\beta-\gamma}
\overline{\widetilde\lambda_{\alpha-\beta-\delta}}=\nonumber
\\
\sum\limits_{\nul\le\gamma\le\alpha}\sum\limits_{\gamma\le\epsilon\le\alpha}
\sum\limits_{\gamma\le\beta\le \alpha-\epsilon+\gamma}
\lll\alpha\atop\beta\rrr\lll\beta\atop\gamma\rrr\lll\alpha-\beta\atop\epsilon-\gamma\rrr
(-\odin)^\gamma(i\rho)^\epsilon\lambda_{\beta-\gamma}
\overline{\widetilde\lambda_{\alpha-\beta-\epsilon+\gamma}}=\nonumber
\\
\sum\limits_{\nul\le\epsilon\le\alpha}(i\rho)^\epsilon
\sum\limits_{\nul\le\gamma\le\epsilon}
\sum\limits_{\gamma\le\beta\le \alpha-\epsilon+\gamma}
\lll\alpha\atop\beta\rrr\lll\beta\atop\gamma\rrr\lll\alpha-\beta\atop\epsilon-\gamma\rrr
(-\odin)^\gamma\lambda_{\beta-\gamma}
\overline{\widetilde\lambda_{\alpha-\beta-\epsilon+\gamma}}=\nonumber
\\
\sum\limits_{\nul\le\epsilon\le\alpha}(i\rho)^\epsilon
\sum\limits_{\nul\le\gamma\le\epsilon}
\sum\limits_{\nul\le\kappa\le \alpha-\epsilon}
\lll\alpha\atop\kappa+\gamma\rrr\lll\kappa+\gamma\atop\gamma\rrr\lll\alpha-\kappa-\gamma\atop\epsilon-\gamma\rrr
(-\odin)^\gamma\lambda_{\kappa}
\overline{\widetilde\lambda_{\alpha-\kappa-\epsilon}}=\nonumber
\\
\sum\limits_{\nul\le\epsilon\le\alpha}(i\rho)^\epsilon
\sum\limits_{\nul\le\kappa\le \alpha-\epsilon}
\lambda_{\kappa}\overline{\widetilde\lambda_{\alpha-\kappa-\epsilon}}
\sum\limits_{\nul\le\gamma\le\epsilon}
\frac{\alpha!}{\kappa!\gamma!(\epsilon-\gamma)!(\alpha-\kappa-\epsilon)!}
(-\odin)^\gamma=\nonumber
\\
\sum\limits_{\nul\le\epsilon\le\alpha}(i\rho)^\epsilon
\sum\limits_{\nul\le\kappa\le \alpha-\epsilon}
\lambda_{\kappa}\overline{\widetilde\lambda_{\alpha-\kappa-\epsilon}}
\lll\alpha\atop\kappa+\epsilon\rrr\sum\limits_{\nul\le\gamma\le\epsilon}
\lll\epsilon\atop\gamma\rrr(-\odin)^\gamma=\nonumber
\\
\sum\limits_{\nul\le\epsilon\le\alpha}(i\rho)^\epsilon
\sum\limits_{\nul\le\kappa\le \alpha-\epsilon}
\lambda_{\kappa}\overline{\widetilde\lambda_{\alpha-\kappa-\epsilon}}
\lll\alpha\atop\kappa+\epsilon\rrr(\odin-\odin)^\epsilon=
\sum\limits_{\nul\le\kappa\le \alpha}\lll\alpha\atop\kappa\rrr
\lambda_{\kappa}\overline{\widetilde\lambda_{\alpha-\kappa}}.\
\label{63}
\ea
So, $D^\alpha\lll \mu_{0 k}(x)\overline{\widetilde\mu_{0 k}(x)} \rrr\Big|_{x=\nul}$
does not depend on $k$ and
\ban
\sum\limits_{\nul\le\kappa\le \alpha}\lll\alpha\atop\kappa\rrr\lambda_{\kappa}\overline{\widetilde\lambda_{\alpha-\kappa}}=
\sml k0{m-1}D^\alpha\lll \mu_{0 k}(x)\overline{\widetilde\mu_{0 k}(x)} \rrr\Bigg|_{x=\nul}=
D^\alpha\lll \sml k0{m-1}\mu_{0 k}(x)\overline{\widetilde\mu_{0 k}(x)} \rrr\Bigg|_{x=\nul}.
\ean

Due to Lemmas~\ref{l1}, \ref{l2}, the $(r+1)\times m$ matrixes ${\cal M}, \widetilde{\cal M}$
can be extended to $(r+1)\times (r+1)$ matrixes ${\cal N}, \widetilde{\cal N}$ such that their entries
are in $L^{(n)}_\infty$ and
${\cal N}^T\overline{\widetilde{\cal N}}=I_{r+1}$.
Due to  conditions (b) of Theorem~\ref{t2},
$$
D^\alpha\lll \mu_{0 k}(x)\overline{\widetilde\mu_{0 k}(x)} \rrr\Bigg|_{x=\nul}=0, \ k=m\ddd r.
$$
From this, taking into account that $\sml k0{r}\mu_{0 k}(x)\overline{\widetilde\mu_{0 k}(x)}=1$,
we obtain
\ban
D^\alpha\lll \sml k0{m-1}\mu_{0 k}(x)\overline{\widetilde\mu_{0 k}(x)} \rrr\Bigg|_{x=\nul}=
D^\alpha\lll \sml k0{r}\mu_{0 k}(x)\overline{\widetilde\mu_{0 k}(x)} \rrr\Bigg|_{x=\nul}=0. \Diamond
\ean

\begin{coro}
Let  $\mu_{0k},\widetilde\mu_{0k}\in L_\infty^{(n)}$ for some $k=0\ddd m-1$.
If  there exist complex numbers $\lambda_\gamma, \widetilde\lambda_\gamma$,
$\gamma\in\zd_+$,  $[\gamma]\le n$, such that~(\ref{58}), (\ref{59}), (\ref{57}) are fulfilled,
 then
$$D^\alpha\lll\overline{\mu_{0k}(x)}\widetilde\mu_{0k}(x)\rrr\Big|_{x=\nul}=0,\ \ \
$$  for all
  $\alpha\in \zd_+$, $0<[\alpha]\le n$.
  \label{c1}
\end{coro}

The proof of Corollary~\ref{c1}  follows from~(\ref{63}).

\subsection {4. Methods for construction  compactly supported
wavelet frames with $VM^n$ property.}

Due to Corollary~\ref{c1} and Theorem~\ref{t5},  refinable functions
$\phi,  \widetilde\phi$ whose masks are trigonometric polynomials
satisfying~(\ref{58}), (\ref{59}), (\ref{57}) generate dual
compactly supported  wavelet systems $\{\psi_{jk}^{(\nu)}\}$,
$\{\widetilde\psi_{jk}^{(\nu)}\}$ with  $VM^n$ property for
$\{\psi_{jk}^{(\nu)}\}$ (we can set
$\mu_{0m}\equiv1, \widetilde\mu_{0m}=\sml l0{m-1}\overline{\mu^\prime_{0l}}\,\widetilde\mu_{0l}$,
in this case). But this is a bad construction. The system $\{\widetilde\psi_{jk}^{(\nu)}\}$
does not have  $VM^0$ property. As was mentioned above, such a system can not be a frame.
The following method allows to provide $VM^n$ property for each of the systems $\{\psi_{jk}^{(\nu)}\}$,
$\{\widetilde\psi_{jk}^{(\nu)}\}$.

{\bf Step 1.} Given $n\in\zd$ and given a set of parameters $\lambda_\beta\in{\Bbb C}$,
$\beta\in \zd_+$, $[\beta]\le n$, $\lambda_\nul=1$,
 find a dual set of parameters $\widetilde\lambda_\beta\in{\Bbb C}$ satisfying~(\ref{57})
 by the following recursive  formulas
$$
\widetilde\lambda_\nul=1,\ \ \ \ {\widetilde\lambda_\alpha}=-\overline{\lambda_\alpha}
-\sum\limits_{\nul<\beta\le\alpha}
\lll\alpha\atop\beta\rrr\overline{\lambda_\beta}{\widetilde\lambda_{\alpha-\beta}}.
$$

{\bf Step 2.} Chose functions $\mu^\prime_{00}\ddd \mu^\prime_{0\,m-1}$ and
$\widetilde\mu_{00}\ddd \widetilde\mu_{0\,m-1}$ defined by
\ban
\mu^\prime_{0k}(x)=
\frac1{\sqrt m}\sum\limits_{[\alpha]\le n}g_\alpha(x)
\sum\limits_{\nul\le\beta\le\alpha}
\lll\alpha\atop\beta\rrr\lll-2\pi i M^{-1}s_k)\rrr^{\beta}\lambda_{\alpha-\beta}+
\\
\sum\limits_{[\alpha]= n+1}T_{k,\alpha}(x)\prod\limits_{j=1}^{\alpha_j}\lll 1-\ex{x_j}\rrr^{\alpha_j},
\ean
\ban
\widetilde\mu_{0k}(x)=
\frac1{\sqrt m}\sum\limits_{[\alpha]\le n}g_\alpha(x)
\sum\limits_{\nul\le\beta\le\alpha}
\lll\alpha\atop\beta\rrr\lll-2\pi i M^{-1}s_k)\rrr^{\beta}\widetilde\lambda_{\alpha-\beta}+
\\
\sum\limits_{[\alpha]= n+1}\widetilde T_{k,\alpha}(x)\prod\limits_{j=1}^{\alpha_j}\lll 1-\ex{x_j}\rrr^{\alpha_j},
\ean
where $T_{k,\alpha}, \widetilde T_{k,\alpha}$ are arbitrary trigonometric polynomials,
$g_\alpha$ are trigonometric polynomials such that
$D^\alpha g_\alpha(\nul)=1$, $D^\beta g_\alpha(\nul)=0$
for all $\beta\in\zd_+$, $\beta\ne\alpha$, $[\beta]\le n$ (recursive formulas for computing $g_\alpha$
are given in~\cite{36}). It is clear that~(\ref{59}) are fulfilled

{\bf Step 3.}
Set $\sigma:=\sml l0{m-1}\overline{\mu^\prime_{0l}}\,\widetilde\mu_{0l}$,
\ \ \ $\mu_{0k}:=(2-\sigma)\mu^\prime_{0k}, k=0\ddd m-1$.

Due to Corollary~\ref{c1}, we have $D^\beta\sigma(\nul)=0$
for all   $\beta\in\zd_+$,  $0<[\beta]\le n$. It follows that
$D^\beta\mu_{0k}(\nul)=D^\beta\mu^\prime_{0k}(\nul)$
for all  $\beta\in\zd_+$,  $[\beta]\le n$.
It is not difficult to see that~(\ref{58}) holds and
$$
1-\sml k0{m-1}\overline{\mu_{0k}(x)}\widetilde\mu_{0k}=(1-\sigma)^2.
$$

Set $\mu_{0m}:=1-\sigma, \widetilde\mu_{0m}:=\overline{1-\sigma}$.

{\bf Step 4.}\
Find  matrixes
$$
\cal M=\lll\matrix{
\mu_{00}&\dots&\mu_{0, m-1}&\mu_{0, m}
\cr
\mu_{10}&\dots&\mu_{1, m-1}&*
\cr
\dots&\dots&\dots&\dots
\cr
\mu_{m,0}&\dots&\mu_{m, m-1}&*
}\rrr,\ \ \
\widetilde {\cal M}=\lll\matrix{
\widetilde \mu_{00}&\dots&\widetilde \mu_{0, m-1}&\widetilde \mu_{0, m}
\cr
\widetilde\mu_{10}&\dots&\widetilde\mu_{1, m-1}&*
\cr
\dots&\dots&\dots&\dots
\cr
\widetilde \mu_{m,0}&\dots&\widetilde\mu_{m, m-1}&*
}\rrr
$$
such that their entries are trigonometric polynomials and  ${\cal M}\overline{\widetilde {\cal M}^T}=I_{m+1}$.

Though the matrixes $\cal M, \widetilde {\cal M}$ can be constructed theoretically (see the proof of
Theorem~\ref{t5}), it is very complicate to implement the algorithm in practice.
Instead, we suggest the following explicit way (the payment of  simplicity of this way is
increasing of the redundancy).

Set $\mu_{0,m+1}\equiv0, \widetilde\mu_{0, m+1}\equiv0$. For each $\nu=1\ddd m+1$, define
\ban
\widetilde\mu_{\nu, m+1}:=\overline{\mu_{0,  m+1-\nu}},\ \
\mu_{\nu, m+1}:=\overline{\widetilde\mu_{0,  m+1-\nu}},\ \ \hspace{6.7cm}
\\
\mu_{\nu k}:=\delta_{m+1-\nu, k}- {\mu_{0 k
}\overline{\widetilde\mu_{0, m+1-\nu}}}, \ \ \ \widetilde\mu_{\nu
k}:=\delta_{m+1-\nu, k}- {\widetilde\mu_{0 k }\overline{\mu_{0,
m+1-\nu}}}, \ \ \    k=0\ddd m. \ean It is not difficult to see
that the matrixes ${\cal M}:=\{\mu_{\nu k}\}_{\nu,k=0}^{m+1}$,
$\widetilde{\cal M}:=\{{\widetilde \mu_{\nu k}}\}_{\nu,k=0}^{m+1}$
satisfy ${\cal M}\overline{\widetilde {\cal M}^T}=I_{m+2}$.

Realization of suggested construction is simple, and there is a good chance to succeed
in applications using the corresponding  dual wavelet systems $\psi^{(\nu)}_{jk},
\widetilde \psi^{(\nu)}_{jk}$. However, as it was mentioned above, we can  be sure that
these systems form dual frames only if additional conditions are fulfilled, in particularly, if
 $\h\phi, \h{\widetilde\phi}$ decay fast enough. Checking of such
 conditions is very complicate in practice. Moreover, usually it is also difficult to check
that $\phi, {\widetilde\phi}\in\Ldvad$. Situation is  better for the case $\psi^{(\nu)}_{jk}=
\widetilde \psi^{(\nu)}_{jk}$. Namely, if $\phi\in\Ldvad$ is a refinable function with a mask $m_0$,
$\lim\limits_{\xi\to0}\h\phi(\xi)=1$, periodic functions  $m_\nu$,
$\nu=1\ddd r$,  are so that the corresponding polyphase matrix
\ban
{\cal M}:=\left(%
\begin{array}{ccc}
  \mu_{00} & \dots& \mu_{0, m-1} \\
  \vdots & \ddots & \vdots\\
  \mu_{r,0} & \dots & \mu_{r,m-1}\\
\end{array}%
\right)
\ean
satisfies
\be
{\cal M^T}\overline{{\cal M}}=I_m,
\label{74}
\ee
and 
$
\widehat\psi^{(\nu)}(x)=m_\nu({M^*}^{-1}x) \widehat\phi({M^*}^{-1}x),
$
then  $\{\psi_{jk}^{(\nu)}\}$ is a tight frame. This fact was presented in~\cite[Corollary 6.7]{RS}
with an additional assumption  which may be omitted if we improve
the proof repeating one-dimensional arguments of~\cite{10}.
So, for any appropriate  polyphase matrix ${\cal M}$, we may be sure that
$\{\psi_{jk}^{(\nu)}\}$ is a tight frame whenever
 $\phi\in\Ldvad$ and $\h\phi$ is continuous at $\nul$.
 If $m_0$ is a trigonometric polynomial, then $\phi$ is compactly supported, which yields continuity
 of $\h\phi$   whenever  $\phi\in\Ldvad$. But this condition will be also fulfilled
 because, due to~(\ref{74}), we have
 \be
\sml k0{m-1}|\mu_{0 k}|^2\le1,
 \label{70}
 \ee
 and the following generalization of well known Mallat's theorem holds.

\begin{prop} Let the polyphase representatives of a trigonometric polynomial~$m_0$
satisfy~(\ref{70}), $\mu_{0k}(\nul)=1/\sqrt m$, $k=0\ddd m-1$,
and let
$f(x):=\prod\limits_{j=1}^{\infty}{m_0({M^*}^{-j}x)}$.
Then  $f\in\Ldvad$, $\|f\|\le1$.
\label{1002}
\end{prop}

{\tt Proof.} First of all we note  (see, e.g.,~\cite{NPS}, \S2.6) that
\be
\sum \limits_{s \in D(M^*)}| m_0(x+ {M^*}^{-1}s)|^2  =
\sum \limits_{k=0}^{m-1}|\mu_{0 k}({M^*}x)|^2.
\label{d72}
\ee
\label{dl6}
It follows from~(\ref{70}) that
\be
\sum \limits_{s \in D(M^*)}| m_0(x+ {M^*}^{-1}s)|^2\le1.
\label{72}
\ee
Set
$
f_k(x)= \prod\limits_{j=1}^{k}{m_0({M^*}^{-j}x)}\chi_{{M^*}^{k}[-1/2,1/2]^d}(x).
$

 Taking into account  that for any $1$-periodic (on each variable) function $g$
$$
\int\limits_{[0,1]^d}g(x)\,dx=\sum \limits_{r\in
D({M^*})}\int\limits_{{M^*}^{-1}[0,1)^d+{M^*}^{-1}r}g(x)\,dx
$$
 (see, e.g., \cite[\S~2.2]{NPS} ), we have
\ban \|f_k\|^2=
\int\limits_{{M^*}^k[-1/2,1/2]^d} \prod
\limits_{j=1}^{k}|{m_0({M^*}^{-j}x)}|^2\,dx=
m^k\int\limits_{[0,1]^d} \prod
\limits_{j=1}^{k}|{m_0({M^*}^{k-j}x)}|^2\,dx=
\\
m^k\sum \limits_{r\in
D({M^*})}\int\limits_{{M^*}^{-1}[0,1)^d+{M^*}^{-1}r} \prod
\limits_{j=1}^{k}|{m_0({M^*}^{k-j}x)}|^2\,dx=
\\
m^k\sum \limits_{r\in D({M^*})}\int\limits_{{M^*}^{-1}[0,1)^d}
\prod \limits_{j=1}^{k-1}|{m_0({M^*}^{k-j}x)}|^2|{m_0(x+
{M^*}^{-1} r)}|^2\,dx\le
\\
m^k\int\limits_{{M^*}^{-1}[0,1)^d} \prod
\limits_{j=1}^{k-1}|{m_0({M^*}^{k-j}x)}|^2\,dx=
m^{k-1}\int\limits_{[0,1)^d} \prod
\limits_{j=1}^{k-1}|{m_0({M^*}^{k-1-j}x)}|^2\,dx=
\\
\int\limits_{{M^*}^{k-1}[0-1/2,1/2]^d} \prod
\limits_{j=1}^{k-1}|{m_0({M^*}^{-j}x)}|^2\,dx
=\|f_{k-1}\|^2,
\ean
Since, by~(\ref2), $m_0(\nul)= 1$, the infinit product converges at each point,
which yields $f_k(x)\too\limits_{k\to+\infty}f(x)$ for all $x\in\rd$. It follows from Fatou's lemma
$$
\int\limits_{\rd}|f(x)|^2\,d\xi=\int\limits_{\rd}
\prod \limits_{j=1}^{\infty}|{m_0({M^*}^{-j}\xi)}|^2\,d\xi=
\sup_{k\ge0}\int\limits_{\rd}|f_k(\xi)|^2\,d\xi\le1.\Diamond
$$

Now let us discuss how to choose a suitable refinable mask  $m_0$ for constructing compactly supported
tight wavelet frames with $VM^n$ property.  Note one more necessary condition for  $VM^n$ property.
If $\{\psi_{jk}^{(\nu)}\}$ is a tight frame and $\cal M$ is its polyphase matrix with polynomial
entries, then, due to Suslin's theorem\cite{3} and
Lemma~\ref{l2}, there exist $(r+1)\times (r+1)$ polynomial matrixes
${\cal N}=\{\mu_{\nu k}\}_{\nu,k=0}^r, \widetilde{\cal N}=\{\widetilde\mu_{\nu k}\}_{\nu,k=0}^r$
extending $\cal M$ such that
${\cal N}\overline{{\widetilde{\cal N}}^T}=I_{r+1}$.
If  $VM^n$ property holds for  $\{\psi_{jk}^{(\nu)}\}$,
by Theorems~\ref{t2} and \ref{t6},
\be
D^\beta\mu_{0k}(\nul)=\frac1{\sqrt m}
\sum\limits_{\nul\le\gamma\le\beta}\lambda_\gamma
\lll\beta\atop\gamma\rrr(-2\pi i M^{-1}s_k)^{\beta-\gamma} \ \ \
\forall    \beta\in\zd_+,  [\beta]\le n,
\label{73}
\ee
where a set of parameters $\lambda_\alpha$ is so that
\be
\lambda_\nul=1,\ \ \
\sum\limits_{\nul\le\gamma\le \alpha}\lll\alpha\atop\gamma\rrr
\lambda_\gamma \overline{\lambda_{\alpha-\gamma}}=0 \ \ \
\forall \alpha\in\zd_+, 0<[\alpha]\le n.
\label{75}
\ee
Of course, providing (\ref{73}) and (\ref{75}) is not enough, in particular,
these conditions do not guarantee (\ref{70}). Now we will show that starting with
arbitrary trigonometric polynomials satisfying (\ref{73}), (\ref{75}),
it is possible to improve them such that the improved polynomials
are polyphase representatives of a refinable mask generating a
tight frame with $VM^n$ property. Our algorithm is based on
the following result by M.A.Dritschel.
\begin{theo}\cite{50}
If $T$ is a strictly  positive trigonometric polynomial of $d$~variables, $d\ge1$, then
there exist trigonometric polynomials $t_1\ddd t_N$ so that
\be
T=|t_1|^2+\dots+ |t_N|^2.
\label{91}
\ee
\label{t7}
\end{theo}

Note that another proof of this theorem was found by J.~Jeronimo and M.~J.~Lai~\cite{111}.
Their proof is constructive for the case $d=2$.
In the case $d=1$, due to Riesz lemma,  (\ref{91}) holds with $N=1$
for any non-negative trigonometric
polynomial $T$.

{\bf Step 1.} Given  parameters $\lambda_\alpha$,  $ 0<[\alpha]\le n$,  which satisfy~(\ref{75}).
choose functions $\mu^\prime_{0k}$, $k=1\ddd m-1$,  defined by
\ba
\mu^{\prime}_{0 k}(x)=
\frac1{\sqrt m}\sum\limits_{[\alpha]\le 2n}g_\alpha(x)
\sum\limits_{\nul\le\beta\le\alpha}
\lll\alpha\atop\beta\rrr\lll-2\pi i M^{-1}s_k)\rrr^{\beta}\lambda_{\alpha-\beta}+
\nonumber
\\
\sum\limits_{[\alpha]= 2n+1}T_{k, \alpha}(x)\prod\limits_{j=1}^{\alpha_j}\lll 1-\ex{x_j}\rrr^{\alpha_j},\ \ \
k=1\ddd m-1,
\label{79}
\ea
where $T_{k, \alpha}$ are arbitrary trigonometric polynomials,
$g_\alpha$ are trigonometric polynomials such that
$D^\alpha g_\alpha(\nul)=1$, $D^\beta g_\alpha(\nul)=0$
for all $\beta\in\zd_+$, $\beta\ne\alpha$, $[\beta]\le n$ (see~\cite{36}).

{\bf Step 2.}
Choose large enough positive integers $L_i, M_i$, $i=1\ddd d$, such that the functions
$$
\mu^{\prime\prime}_{0 k}(x):=\prod\limits_{i=1}^d\lll1-\sin^{2L_i}
\pi x_i\rrr^{M_i}\mu^\prime_{0k}(x),\  k=0\ddd m-1,
$$
satisfy
 \be
\sml k0{m-1}|\mu^{\prime\prime}_{0 k}|^2\le2 \ \ \
 \label{90}
 \ee
 and $D^\beta\mu^{\prime\prime}_{0k}(\nul)=D^\beta\mu^\prime_{0k}(\nul)$,
$k=1\ddd d$, for all $ \beta\in\zd_+, [\beta]\le n$.

{\bf Step 3.}
Set $\sigma:=\sum\limits_{k=0}^{m-1}|\mu^{\prime\prime}_{0k}|^2$,
$$
\mu_{0 k}(x):=\lll\frac32-\frac\sigma2\rrr\mu^{\prime\prime}_{0k}(x),
\  k=0\ddd m-1,
$$
Due to Corollary~\ref{c1}, $D^\beta\mu_{0k}(\nul)=
D^\beta\mu^{\prime\prime}_{0k}(\nul)=D^\beta\mu^\prime_{0k}(\nul)$,
$k=1\ddd d$, for all $ \beta\in\zd_+, [\beta]\le n$.

{\bf Step 4.}
Find  trigonometric polynomials $\mu_{0m}\ddd \mu_{0r}$  so that
\ba
&&1-\sml k0{m-1}|\mu_{0 k}|^2\equiv\sml km{r}|\mu_{0 k}|^2,
\label{81}
\\
&&D^\beta\mu_{0k}(\nul)=0,\ \ \ k=m\ddd r,\ \ \ \forall    \beta\in\zd_+,  [\beta]\le n.
\label{82}
\ea

Such polynomials exist due to the identity
$$
1-\lll\frac32-\frac\sigma2\rrr^2\sigma
=(1-\sigma)^2\lll1-\frac\sigma4\rrr,
$$
the fact that  $1-\frac\sigma4\ge \frac12$ (because of~(\ref{90})), Theorem~\ref{t7}
and Corollary~\ref{c1}.

{\bf Step 5.}
Find a unitary matrix
$$
\cal M=\lll\matrix{
\mu_{00}&\dots&\mu_{0, m-1}&\mu_{0m}&\dots&\mu_{0r}
\cr
\mu_{10}&\dots&\mu_{1, m-1}&*&\dots&*
\cr
\dots&\dots&\dots&\dots&\dots&\dots
\cr
\mu_{r,0}&\dots&\mu_{r, m-1}&*&\dots&*
}\rrr,\ \ \
$$
whose entries are trigonometric polynomials.

It is not known if any appropriate row may be extended to a unitary matrix with
polynomial entries. But if $\mu_{0r}\equiv0$ (which  always can be realized),
then $\cal M$ can be defined by:  for all
  $\nu=1\ddd r$ set
\be
\mu_{\nu r}:=\overline{\mu_{0,  r-\nu}},\ \
 \ ,\ \ \
\mu_{\nu k}:=\delta_{r-\nu, k}- {\mu_{0 k
}\overline{\mu_{0, r-\nu}}}, \ \ \    k=0\ddd r.
\label{83}
\ee

 \vspace{.5cm}
\noindent
Maria Skopina, Department of Applied Mathematics and Control
Processes, Saint Petersburg State University,

\noindent
e-mail:\ skopina@MS1167.spb.edu

\begin {thebibliography} {99}

\bibitem{9}
A.P.Calder\'on, Zygmund  A., {\em Local properties of
solutions of elliptic partial differential equation,}
Studia Math., 20 (1961), p. 171-227.



\bibitem{15}
Chen Di.-R.,  Han H. and  Riemenschneider S. D.
{\it Construction of multivariate biorthogonal wavelets with arbitrary vanishing moments},
Adv. Comput. Math. {\bf 13 } (2000), No.2, 131-165.


\bibitem{35}
Ingrid Daubechies, Bin Han, Amos Ron, and Zuowei Shen,
{\it Framelets: MRA-based constructions of wavelet frames},
ACHA, Vol. 14 (2003), No.~1, 1--46.

\bibitem{50}
M.A.Dritschel {\it On factorization of trigonometric polynomials},
Integral Equations and Operator Theory, 49(2004), 11--42.


\bibitem{111}
 Geronimo J. and Lai M.J.{\it  Factorization of Multivariate Positive  Laurent Polynomials}
 J. of Approx. Theory 139(2006), 327--345.


\bibitem{98}
Bin Han. {\it On dual wavelet tight frames} ACHA. 1997. V.~4. P.~380-413.

\bibitem {b16}
Hernandes E. and Weis G.A. {\it A first cours of wavalets}
CRC Press, Boca Raton, FL. 1996.

\bibitem{32}
Ji H., Riemenschneider S.D. and Shen Z.
{\it Multivariate compactly supported fundamental
refinable functions, dual and biorthogonal wavelets},
Studies of Applied Mathematics {\bf 102 }(1999), 173-204.

\bibitem{1}
Jia R.Q. {\it Approximation properties of multivariate wavelets},
Math. Comp. {\bf  67 }(1998), 647-655.


\bibitem {MP}
Ming-Jun Lai and  Petukhov A. {\it Method of Virtual Components
for Constructing Redundant Filter and Wavelet Frames} (to appear).

\bibitem {11}
M\"oller H.M. and  Sauer T.
{\it Multivariate refinable functions of high approximation
order via quotient ideals of Laurent polynomials},
Adv. Comput. Math. {\bf 20 } (2004),  No.1-3, 205-228.


\bibitem{NPS}
Novikov I., Protassov V., and Skopina M.  {\it Wavelet Theory}. Moscow: Fizmatlit, 2005.


\bibitem{10}
Petukhov A. {\it Explicit construction of framelets}
 ACHA, Vol.~11 (2001(, ~313-327.


\bibitem{P}
Petukhov   A. {\it Construction of Symmetric Orthogonal Bases of Wavelets and
          Tight Wavelet Frames  with Integer Dilation Factor},
          Appl. Comput. Harmon. Anal., {\bf 17} (2004), 198--210.

\bibitem{40}
Riemenschneider S.D. and Shen Z.W. {\it Construction of compactly
supported biorthogonal wavelets in $L_2 (\r^s) $}, Preprint, 1997.



\bibitem{RS}
Ron  A. and  Shen Z., {\it Affine systems in ${\mathbb L_2(R^d)}$:
the analysis of the analysis operator}, J. Func. Anal.  {\bf 148}
(1997), 408--447.




\bibitem{36}
Skopina M.{\it On Construction of Multivariate Wavelets with Vanishing Moments},
ACHA, Vol. 20 (2006),  3  375-390.


\bibitem{3}
Suslin A. {\it The structure of the special linear group over rings of polynomials},
Izv. Akad. Nauk SSSR. Ser. Mat. {\bf 41} (1977), No 2,  235-252 (in Russian)




 \bibitem {80}
Shen Z. Extension of matrices with Laurent polynomial entries~//
Proceedings of the 15th IMACS World Congress on Scientific Computation
  Modeling and Applied Mathematics, Ashim Syclow eds. 1997. P.~ 57-61.




\end {thebibliography}
\end{document}